\documentclass[12pt]{amsart}
\usepackage[utf8]{inputenc}
\usepackage{graphicx} 

\usepackage{graphicx}
\usepackage{amssymb, amsthm, amsmath, tcolorbox}
\usepackage{enumitem}
\usepackage{tikz-cd}
\usepackage{bbold}
\usepackage{ytableau}   
\usepackage[margin=1in]{geometry}
\usepackage{ifthen}

\usepackage{hyperref}

\theoremstyle{definition}
\newtheorem*{theorem*}{Theorem}

\newtheorem{theorem}{Theorem}[section]
\newtheorem{lemma}[theorem]{Lemma}
\newtheorem{proposition}[theorem]{Proposition}

\newtheorem{remark}[theorem]{Remark}
\newtheorem{example}[theorem]{Example}
\newtheorem{conjecture}[theorem]{Conjecture}

\newtheorem{corollary}[theorem]{Corollary}
\newtheorem{definition}[theorem]{Definition}

\newcommand{\BK}{\operatorname{BK}}
\newcommand{\ssl}{\mathfrak{sl}}

\def\Imm{\operatorname{Imm}}
\def\TL{\operatorname{TL}}

\newtheorem{innercustomthm}{Theorem}
\newenvironment{repthm}[1]
  {\renewcommand\theinnercustomthm{#1}\innercustomthm}
  {\endinnercustomthm}

\newcommand{\uncross}[3]{
    \draw[thick] [bend right = 90, looseness=1.25] (-0.5+#1, -#2) to (-0.5+#1, -#2+1);
    \draw[thick] [bend right = 90, looseness=1.25] (0.5+#1, -#2+1) to (0.5+#1, -#2);
    \foreach \j in {0,...,#3}{
    \ifthenelse{\equal{\j}{\the\numexpr #2-1} \OR \j = #2}{}{
    \draw[thick] (-0.5+#1,-\j ) to (0.5+#1,-\j);}
    }
}

\title{Shuffle Tableaux, Littlewood--Richardson Coefficients, and Schur Log-Concavity}
\author{Chau Nguyen, Son Nguyen, Dora Woodruff}

\begin{document}

    \begin{abstract}
        We give a new formula for the Littlewood--Richardson coefficients in terms of peelable tableaux compatible with shuffle tableaux, in the same fashion as Remmel--Whitney rule. This gives an efficient way to compute generalized Littlewood--Richardson coefficients for Temperley--Lieb immanants of Jacobi--Trudi matrices. We will also show that our rule behaves well with Bender--Knuth involutions, recovering the symmetry of Littlewood--Richardson coefficients. As an application, we use our rule to prove a special case of a Schur log-concavity conjecture by Lam--Postnikov--Pylyavskyy.
    \end{abstract}

\maketitle

\tableofcontents

\section{Introduction}

    The Littlewood--Richardson coefficients $c_{\lambda,\nu}^\nu$ are some of the most important structure constants in algebraic combinatorics. In combinatorics, they are the structure constants for the product of Schur functions $s_\lambda s_\mu = \sum_\nu c_{\lambda,\mu}^\nu s_\nu$. In representation theory, $c_{\lambda,\mu}^{\nu}$ is the multiplicity of the irreducible representation $V_{\nu}$ of $S_{|\nu|}$ in $\operatorname{Ind}_{S_{|\lambda|}\times S_{|\mu|}}^{S_{|\nu|}}V_\lambda\otimes V_\mu$, as well as the multiplicity of the Schur module $E^{\nu}$ in $E^{\lambda}\otimes E^{\mu}$. In Schubert calculus, they are the Schubert structure constants for the product of Schubert classes corresponding to Grassmannian permutations $[X_{\lambda}]\cdot[X_\mu] = \sum_\nu c_{\lambda,\mu}^\nu[X_\nu]$.

    Since the Littlewood--Richardson rule was introduced \cite{littlewood1934group}, several variants of the rule have been formulated, including Knutson--Tao honeycombs \cite{knutson1999honeycomb}, Zelevinsky pictures \cite{zelevinsky1981generalization}, Berenstein--Zelevinsky triangles \cite{berenstein1992triple} and their more general rule for any Lie algebra $\mathfrak{g}$ \cite{berenstein1999tensor}. A somewhat similar rule to Zelevinsky pictures is Remmel--Whitney rule \cite{remmel1984multiplying}, which can be formulated in terms of counting Reiner--Shimozono \textit{peelable tableaux} \cite{reiner1998percentage} as follows. Let $D$ be the skew diagram of shape $\lambda\ast\mu = (\lambda_1+\mu_1,\ldots,\lambda_k +\mu_1,\mu_1,\ldots,\mu_\ell)/(\mu_1^k)$. Visually, $D$ is obtained by putting the Young diagram of $\lambda$ to the top right of that of $\mu$. For $1\leq i\leq \ell(\nu)-1$, let $a_i$ be the number of columns with squares on both row $i$ and $i+1$. A semistandard Young tableau (SSYT) $T$ is a $D$-peelable tableau if for all $1\leq i\leq \ell(\nu)-1$, there are at least $a_i$ disjoint pairs of $i$-square and $(i+1)$-square such that the $i$-square is Northeast of the $(i+1)$-square. Then, $c_{\lambda,\mu}^{\nu}$ is the number of $D$-peelable tableaux of shape $\nu$.
    
    In this paper, we give a new formula for Littlewood--Richardson coefficients in terms of peelable tableaux. The key difference is that our rule uses the \textit{shuffle diagram} of $\lambda$ and $\mu$ instead of the skew diagram $\lambda * \mu$. In fact, our formula can be used to compute the Schur expansion of any product of two skew-Schur functions. Since shuffle tableaux are used to give a generalized Littlewood--Richardson rule for Temperley--Lieb immanants of Jacobi--Trudi matrices, our rule gives a more efficient way to compute these coefficients. Furthermore, as we will see in Section \ref{subsec:intro-log}, our new rule can be used to prove a special case of Lam--Postnikov--Pylyavskyy's Schur log-concavity conjecture.

\subsection{Shuffle diagram, peelable tableaux, and Littlewood--Richardson coefficients}

    \begin{definition}[$D$-peelable tableaux]\label{def:peelable}
        Given two skew shape $\lambda/\mu$ and $\nu/\rho$, we construct the shuffle diagram $D$ of shape $(\lambda/\mu)\circledast(\nu/\rho)$ by interlacing the squares of $\lambda/\mu$ and $\nu/\rho$ (see Definition \ref{def:shuffle} and Example \ref{ex:shuffleDiagram}). Let $a_i$ be the number of columns with squares on both row $i$ and $i+2$. A SSYT $T$ is a $D$-peelable tableau if for all $i$, there are at least $a_i$ disjoint pairs of one $i$-square and one $(i+2)$-square such that the $i$-square is Northeast of the $(i+2)$-square. We call such pair a \textit{NE matching of $i$ and $i+2$}.
    \end{definition}

    \begin{theorem}\label{thm:main-thm}
        Let $D$ be the shuffle diagram of shape $(\lambda/\mu)\circledast(\nu/\rho)$. For any $\kappa$, the Littlewood--Richardson coefficient $c^{\kappa}_{\lambda/\mu, \nu/\rho}$ counts the number of $D$-peelable tableaux of shape $\kappa$.
    \end{theorem}

    \begin{example}\label{ex:shuffleDiagram}

        Let $\lambda/\mu = (3,2)/(1,0)$ and $(\nu/\rho) = (3,3)/(1,0)$. Then the shuffle diagram of shape $(\lambda/\mu)\circledast(\nu/\rho)$ is the right diagram in Figure \ref{fig:shuffleDiagram}.
        
        \begin{figure}[h!]
            \centering
            \includegraphics[scale = 0.8]{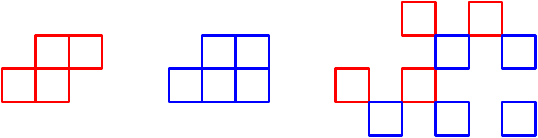}
            \caption{The Young diagram of $(3,2)/(1,0)$ (left) and $(3,3)/(1,0)$ (center) and their shuffle diagram (right)}
            \label{fig:shuffleDiagram}
        \end{figure}

        Thus, for a SSYT $T$ to be $D$-peelable, we need one compatible pair of $1$ and $3$, and two compatible pairs of $2$ and $4$. Let $\kappa = (4,3,2)$, Figure \ref{fig:peelableTableaux} shows all three possible $D$-peelable tableaux of shape $\kappa$. Thus, $c^{\kappa}_{\lambda/\mu, \nu/\rho} = 3$.
        
        \begin{figure}[h!]
            \centering
            \includegraphics[scale = 0.8]{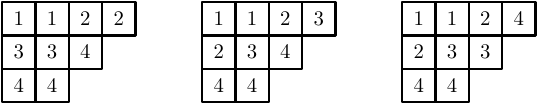}
            \caption{$D$-peelable tableaux}
            \label{fig:peelableTableaux}
        \end{figure}
    \end{example}

\subsection{Temperley--Lieb immanants}

    In his groundbreaking 1990 paper \cite{lusztig1990canonical} Lusztig has introduced {\it {canonical bases}} of quantum groups. They inspired several major developments in mathematics, such as Kashiwara crystal bases \cite{kashiwara1991crystal, kashiwara1993global} and Fomin--Zelevinsky cluster algebras \cite{fomin2002cluster, berenstein1996parametrizations}. In type A, explicit formulas for dual canonical bases were given by Du \cite{du1992canonical, du1995canonical} and by Rhoades--Skandera \cite{rhoades2006kazhdan, skandera2008dual}, where the coefficients are given as evaluations of Kazhdan--Lusztig polynomials. Lusztig observed that they have a deep and surprising relation to total positivity and extended the notion of total positivity to any algebraic group and associated flag varieties \cite{lusztig1994total, lusztig1997total}.

    Another deep positivity property of dual canonical bases in type $A$ is Schur positivity. Any determinant of a (generalized) Jacobi--Trudi matrix is a skew Schur function and hence is Schur positive. Determinants are elements of the dual canonical bases, and it turns out that all of their elements evaluate to Schur positive results on generalized Jacobi--Trudi matrices. This was proved by Rhoades--Skandera in \cite[Proposition 3]{rhoades2006kazhdan}, relying on Haiman's result \cite[Theorem 1.5]{haiman1993hecke}, which in turn relies on the (proofs of) Kazhdan--Lusztig conjecture by Beilinson--Bernstein \cite{BB} and Brylinski--Kashiwara \cite{brylinski1981kazhdan}. However, no direct combinatorial proof was known, and in particular, there is no general interpretation for the coefficients of Schur functions appearing in the evaluation, i.e. the generalized Littlewood--Richardson coefficients.

    Fortunately, for a subset of dual canonical bases, called \textit{Temperley--Lieb immanants}, introduced by Rhoades and Skandera in \cite{rhoades2005temperley}, Nguyen--Pylyavskyy \cite{nguyen2025temperley} gave a combinatorial interpretation in terms of Yamanouchi shuffle tableaux (see Section \ref{subsec:TL-imm}). These Temperley--Lieb immanants have been used by Lam--Postnikov--Pylyavskyy \cite{lam2007schur} to prove several deep Schur-positivity conjectures.

    However, the interpretation given in \cite{nguyen2025temperley} involves checking whether a given shuffle tableaux is Yamanouchi by testing that no lowering crystal operator can be applied. This makes generating Yamanouchi shuffle tableaux relatively inefficient. Using Theorem \ref{thm:main-thm}, Corollary \ref{cor:new-LR-rule} provides a computationally more efficient method to compute the generalized Littlewood--Richardson rule since generating peelable tableaux is generally faster.

\subsection{Bender--Knuth involutions and symmetry of Littlewood--Richardson coefficients}

    The Littlewood--Richardson coefficients have an obvious symmetry $c^{\kappa}_{\lambda/\mu, \nu/\rho} = c^{\kappa}_{\nu/\rho, \lambda/\mu}$. Thus, one may ask for a shape-preserving bijection between $(\lambda/\mu)\circledast(\nu/\rho)$-peelable tableaux and $(\nu/\rho)\circledast(\lambda/\mu)$-peelable tableaux. Surprisingly, the answer is simply Bender--Knuth involutions. Assume for simplicity that $\ell(\lambda) = \ell(\nu)$ (if $\ell(\lambda) > \ell(\mu)$, we can append $0$'s at the end of $\mu$ to get $\ell(\lambda) = \ell(\nu)$), we have the following theorem.

    \begin{theorem}\label{thm:BK-symm}
        Let $T$ be a $(\lambda/\mu)\circledast(\nu/\rho)$-peelable tableau, and let $\ell = \ell(\lambda)$. Let
        \[ T' = \BK_1\circ\BK_3\circ\cdots\circ\BK_{2\ell - 1} (T), \]
        then $T'$ is a $(\nu/\rho)\circledast(\lambda/\mu)$-peelable tableau. Note that $T$ and $T'$ have the same shape.
    \end{theorem}

    \begin{example}\label{ex:BK-symm}

        Applying $\BK_1\circ \BK_3$ to the peelable tableaux in Figure \ref{fig:peelableTableaux}, we get the tableaux in Figure \ref{fig:peelableTableaux2}. Every tableaux in Figure \ref{fig:peelableTableaux2} has at least two compatible pairs of $1$ and $3$ and one compatible pair of $2$ and $4$. Hence, they are all $(\nu/\rho)\circledast(\lambda/\mu)$-peelable, for $\lambda,\mu,\nu,\rho$ defined in Example \ref{ex:shuffleDiagram}.
        
        \begin{figure}[h!]
            \centering
            \includegraphics[scale = 0.8]{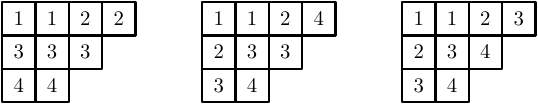}
            \caption{Image under $\BK_1\circ\BK_3$ of the peelable tableaux in Figure \ref{fig:peelableTableaux}}
            \label{fig:peelableTableaux2}
        \end{figure}
    \end{example}

\subsection{Schur log-concavity}\label{subsec:intro-log}

    Recall that the dominant weights of $\ssl_n$ can be viewed as partitions $\lambda$ with $\lambda_n = 0$. Let $V_\lambda$ be the highest weight $\ssl_n$-module corresponding to $\lambda$. A natural question is when a tensor product $V_{\lambda}\otimes V_\mu$ is contained in $V_{\nu}\otimes V_\rho$. Since the multiplicity of $V_{\kappa}$ in $V_{\lambda}\otimes V_\mu$ is exactly $c_{\lambda,\mu}^\kappa$, this question can be phrased as when one has $c_{\lambda,\mu}^\kappa < c_{\nu,\rho}^\kappa$ for all $\kappa$. Combinatorially, this is equivalent to when the difference $s_{\nu}s_\rho - s_{\lambda}s_\mu$ is Schur positive. Several  conjectures and results of this form can be found in \cite{bergeron2004some,fomin2005eigenvalues,lascoux1997ribbon, okounkov1997log}.

    In \cite{lam2007schur}, Lam--Postnikov--Pylyavskyy proved several Schur positivity results, and they conjectured a sufficient condition for this question. Let $\alpha_{ij} = e_i-e_j$ be the roots of the type A root system. An \textit{alcoved polytope} is the intersection of some half-spaces bounded by the hyperplanes $\{\tau~|~\langle \alpha_{ij},\tau\rangle = m\}$. See \cite{lam2007alcoved,lam2018alcoved} for a study of alcoved polytopes. Given two weights $\lambda,\mu$, the \textit{parallelepiped} $P_{\lambda,\mu}$ is the minimal alcoved polytope containing $\lambda$ and $\mu$. We have the following conjecture.

    \begin{conjecture}[{\cite[Conjecture 15.1]{pylyavskyy2007comparing}}]\label{conj:log-concavity}
        If $\lambda + \mu = \nu + \rho$ and $\nu,\rho\in P_{\lambda,\mu}$, then $s_{\nu}s_\rho - s_{\lambda}s_\mu$ is Schur positive.
    \end{conjecture}

    A weaker version of Conjecture \ref{conj:log-concavity} was proved in \cite{dobrovolska2007products}, where they showed that the Schur support of $s_{\lambda}s_\mu$ is contained in that of $s_{\nu}s_\rho$. Here, we resolved a special case of this conjecture by constructing an injection from $\lambda\circledast\mu$-peelable tableaux to $\nu\circledast\rho$-peelable tableaux (Definition \ref{def:theta}). Our injection is relatively simple since it amounts to locally changing $2i-1$ to $2i$ for some $i$. This illustrates the potential of our new rule in proving this conjecture.

    \begin{theorem}\label{thm:log-concavity}
        Suppose $\lambda = (a^k,1^{n-k-1},0)$, and $\nu = \lambda - (e_m + e_{m+1}+\ldots+e_k)$, and $\lambda,\mu,\nu,\rho$ satisfy the conditions of Conjecture \ref{conj:log-concavity}. Then, $s_{\nu}s_\rho - s_{\lambda}s_\mu$ is Schur positive.
    \end{theorem}

    \begin{example}\label{ex:strip}
        Let $\lambda = (7,7,7,7,1,1,1,0)$, $\mu = (9,6,6,5,4,3,3,0)$, $\nu = \lambda - (e_2+e_3+e_4) = (7,6,6,6,1,1,1,0)$, and $\rho = (9,7,7,6,4,3,3,0)$. Then $\lambda,\mu,\nu,\rho$ satisfy the conditions of Theorem \ref{thm:log-concavity}, so $s_{\nu}s_\rho - s_{\lambda}s_\mu$ is Schur positive. Visually, $\nu$ is obtained from $\lambda$ by cutting out a vertical strip from row $2$ to $4$ (see Figure \ref{fig:logConEx}).
        
        \begin{figure}[h!]
            \centering
            \includegraphics[scale = 0.8]{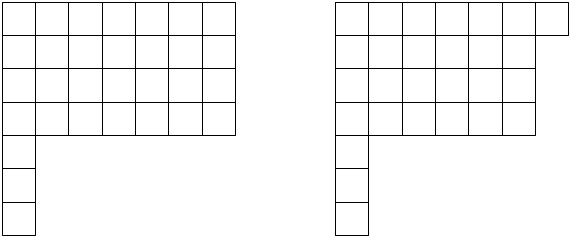}
            \caption{$\lambda$ (left) and $\nu$ (right)}
            \label{fig:logConEx}
        \end{figure}
    \end{example}

    \noindent\textbf{Acknowledgement} We thank Alex Postnikov and Pavlo Pylyavskyy for telling us about their Schur log-concavity conjecture. We thank Vic Reiner for telling us about peelable tableaux. We thank Daniel Soskin for helpful conversations.

\section{Littlewood--Richardson coefficients}

\subsection{Shuffle tableaux and crystal operators}

    First, we review the relevant background on shuffle tableaux and their crystal operators, as introduced in \cite{nguyen2025temperley}.

    \begin{definition}\label{def:shuffle}
        Given two skew shapes $\lambda/\mu$ and $\nu/\rho$, the \textit{shuffle diagram} $D$ of shape $(\lambda/\mu)\circledast(\nu/\rho)$ is constructed as follows.

        \begin{itemize}
            \item If $\lambda/\mu$ has a square in position $(i,j)$, then add into $D$ a square in position $(2i-1,2j-1)$.
            \item If $\nu/\rho$ has a square in position $(i,j)$, then add into $D$ a square in position $(2i,2j)$.
        \end{itemize}
        
        A \textit{shuffle tableau} of shape $(\lambda/\mu)\circledast(\nu/\rho)$ is a filling of the shuffle diagram of shape $(\lambda/\mu)\circledast(\nu/\rho)$ such that the entries are weakly increasing along the rows and strictly increasing along the columns. The \textit{weight} of a shuffle tableau $T$, denoted $\omega(T)$, is $\prod_{i\geq 1} x_i^{c_i}$, where $c_i$ is the number of $i$-entries in $T$.

        In addition, we say square $(r_1,c_1)$ is \textit{right of} square $(r_2,c_2)$ if $r_2 > r_1$ or $r_2 = r_1$ and $c_2 < c_1$. In particular, this means that in the reading order from left to right, bottom to top, square $(r_1,c_1)$ is right of square $(r_2,c_2)$. Then, we also say $(r_2,c_2)$ is \textit{left of} $(r_1,c_1)$.
    \end{definition}

\begin{figure}[h!]
        \centering
        \includegraphics[scale = 0.8]{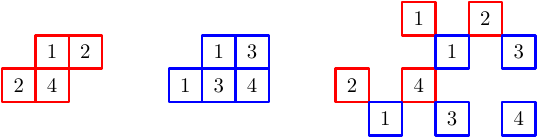}
        \caption{SSYTs of shape $\lambda/\mu$ (left) and $\nu/\rho$ (center) and their corresponding shuffle tableau (right)}
        \label{fig:shuffleTableau}
    \end{figure}

    Observe that the squares in the same row or column of $D$ are exactly those on the same row or column of the Young diagram of either $\lambda/\mu$ or $\nu/\rho$. Thus, one can think of a shuffle tableau of shape $(\lambda/\mu)\circledast(\nu/\rho)$ as interlacing two SSYTs or shape $\lambda/\mu$ and $\nu/\rho$ (see Figure \ref{fig:shuffleTableau}). In particular, we have
    \[ \sum_{\substack{\text{$T$ of shape} \\ (\lambda/\mu)\circledast(\nu/\rho)}} \omega(T) = s_{\lambda/\mu}s_{\nu/\rho}. \]
    
    \begin{definition}
        Given a shuffle tableau $T$, an $(i,i+1)$-overlap is a pair of squares $(s,t)$ on the same column such that $s$ contains an $i$ and $t$ contains an $i+1$. We define the \textit{$i$-reading word} $w_i(T)$ as follows.
        \begin{enumerate}
            \item Consider only the squares containing $i$ and $i+1$.
            \item Remove all $(i,i+1)$-overlap pairs.
            \item Iteratively read the remaining squares from bottom to top, left to right.
        \end{enumerate}
    \end{definition}

    \begin{example}\label{exp:overlap}
        In the following tableau $T$,
        \[ \begin{ytableau}    
        \none & \none & \none & \none & \textcolor{red}{1} & \none & \textcolor{red}{1} & \none & \textcolor{red}{1} & \none & 1 & \none & 2 \\ 
        \none & \none & \none & \textcolor{red}{1} & \none & \textcolor{blue}{2} & \none & 3 & \none & 3 & \none & 3 \\ 
        \none & \none & \textcolor{red}{1} & \none & \textcolor{red}{2} & \none & \textcolor{red}{2} & \none & \textcolor{red}{2} & \none & 3 \\ 
        \none & 2 & \none & \textcolor{red}{2} & \none & \textcolor{blue}{3} \\ 
        2 & \none & \textcolor{red}{2}
        \end{ytableau}, \]
        there are five $(1,2)$-overlap pairs (in red), and one $(2,3)$-overlap pair (in blue). Thus, we have
        \[ w_1(T) = 2~|~2~|~|~2~|~1~2 \]
        by reading all non-red $1$s and $2$s, and
        \[ w_2(T) = 2~2~|~2~2~|~2~2~2~3~|~3~3~3~|~2 \]
        by reading all non-blue $2$s and $3$s.
    \end{example}

    Now we are ready to define the crystal operators on shuffle tableaux.

    \begin{definition}[Crystal operators]\label{def:crys-ops}
        The \textit{crystal operators} $E_i$ and $F_i$ act on $w_i(T)$ in the following (standard) way. 
        \begin{itemize}
            \item View each $i$ as a closing parenthesis ``)'' and each $i+1$ as an opening parenthesis ``(''.
            \item Match the parentheses in the usual way.
            \item $E_i$ changes the leftmost unmatched $i+1$ to $i$ (if exists), and $F_i$ changes the rightmost unmatched $i$ to $i+1$ (if exists).
        \end{itemize}
     
       We induce this action of the crystal operators to the action on shuffle tableaux as follows: 
        \begin{itemize}
            \item the shape of the shuffle tableaux is preserved;
            \item The content changes in the way uniquely determined by $w_i(E_i(T)) = E_i(w_i(T))$ and $w_i(F_i(T)) = F_i(w_i(T))$.
        \end{itemize} 
   \end{definition}

   \begin{example}
       Continuing Example \ref{exp:overlap}, we can view $w_1(T)$ as the following parenthesization.
        \[ \begin{ytableau}
            \none[2] & \none[2] & \none[2] & \none[1] & \none[2]  \\
            \none[(] & \none[(] & \none[\underline{(}] & \none[\underline{)}] & \none[(]
        \end{ytableau} \]
        Hence, $E_1$ changes the first $2$ to $1$, so $E_1(T)$ is the following tableau
        \[ \begin{ytableau}    
        \none & \none & \none & \none & 1 & \none & 1 & \none & 1 & \none & 1 & \none & 2 \\ 
        \none & \none & \none & 1 & \none & 2 & \none & 3 & \none & 3 & \none & 3 \\ 
        \none & \none & 1 & \none & 2 & \none & 2 & \none & 2 & \none & 3 \\ 
        \none & 2 & \none & 2 & \none & 3 \\ 
        \textcolor{red}{1} & \none & 2
        \end{ytableau}, \]
        in which the changed square is colored red.
   \end{example}

   The connected components of the graph on shuffle tableaux formed by the above crystal operators are called \textit{Temperley--Lieb crystals}. In fact, they are type A Kashiwara crystals.

   \begin{theorem}[{\cite[Theorem 5.1]{nguyen2025temperley}}]
       Temperley--Lieb crystals are type A Kashiwara crystals.
   \end{theorem}

   Thus, counting the highest weight tableaux of Temperley--Lieb crystals gives the Littlewood--Richardson coefficients.

   \begin{corollary}\label{cor:yama-shuffle}
       The Littlewood--Richardson coefficient $c_{\lambda/\mu,\nu/\rho}^\kappa$ counts shuffle tableaux of shape $(\lambda/\mu)\circledast(\nu/\rho)$ and content $\kappa$ on which no lowering operator $E_i$ can be applied. These tableaux are called \textit{Yamanouchi shuffle tableaux}.
   \end{corollary}

   \begin{example}\label{ex:yama-shuffle-tableau}
        Recall Example \ref{ex:shuffleDiagram} with $\lambda/\mu = (3,2)/(1,0)$, $(\nu/\rho) = (3,3)/(1,0)$, and $\kappa = 3$. We saw that $c^{\kappa}_{\lambda/\mu, \nu/\rho} = 3$. The three Yamanouchi shuffle tableaux of shape $(\lambda/\mu)\circledast(\nu/\rho)$ and content $\kappa$ are shown in Figure \ref{fig:yamaTableau}.
        
        \begin{figure}[h!]
            \centering
            \includegraphics[scale = 0.8]{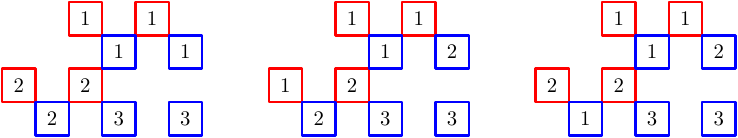}
            \caption{Yamanouchi shuffle tableaux}
            \label{fig:yamaTableau}
        \end{figure}
   \end{example}

   \begin{remark}
       Shuffle tableaux were originally introduced to study Temperley--Lieb immanants, hence the name Temperley--Lieb crystal. In \cite{nguyen2025temperley}, there are additional constraints on the shapes $\lambda/\mu$ and $\nu/\rho$ so that the crystal operators behave nicely with the Temperley--Lieb type (see Section \ref{subsec:TL-imm}). These constraints are unnecessary for this section.
   \end{remark}

\subsection{Compatible tableaux}

    \begin{definition}[$D$-compatible tableaux]\label{def:compatible}
        Let $D$ be the diagram of $(\lambda/\mu)\circledast(\nu/\rho)$. We label the squares of $D$ $1,2,\ldots,|\lambda/\mu| + |\nu/\rho|$ from right to left, bottom to top. We will refer to the square labeled $i$ as square $i$. A SYT $T$ (of any shape) is $D$-compatible if the following two conditions hold:

        \begin{enumerate}
            \item If square $i+1$ is directly left of square $i$ in $D$, then $i+1$ is weakly North and strictly East of $i$ in $T$. 
            \item If square $i$ is directly above square $j$ in $D$, then $i$ is strictly North and weakly East of $j$ in $T$. 
        \end{enumerate}
    \end{definition}

    \begin{example}

        Let $D$ be the same diagram as in Example \ref{ex:yama-shuffle-tableau}. The labeling of $D$ from right to left, bottom to top is shown in Figure \ref{fig:compLabelling}. The two SYTs on the left of Figure \ref{fig:compTableau} are $D$-compatible. The SYT on the right is not because 2 is not weakly North of 1, and 2 is not weakly East of 5.

        \begin{figure}[h!]
            \centering
            \begin{minipage}[b]{.4\textwidth}
                \centering
                \includegraphics[scale = 0.8]{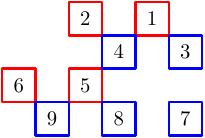}
                \caption{Labeling of $D$}
                \label{fig:compLabelling}
            \end{minipage}%
            \begin{minipage}[b]{.65\textwidth}
                \centering
                \includegraphics[scale = 0.8]{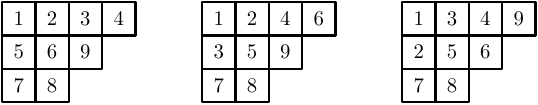}
                \caption{$D$-compatible tableaux (left and center) and non-$D$-compatible tableau (right)}
                \label{fig:compTableau}
            \end{minipage}
        \end{figure}
    \end{example}

    Now, we will construct a bijection $\varphi$ between $D$-compatible tableaux and Yamanouchi shuffle tableaux on $D$. Furthermore, for any $D$-compatible tableau $T$, we would like the shape of $T$ to be the same as the content of $\varphi(T)$.

    \begin{definition}[Construction of $\varphi$]\label{def:yama-comp-bijection}
        Given a $D$-compatible tableau $T$, $\varphi(T)$ is constructed as follows: if $i$ is in row $r$ of $T$, then in $\varphi(T)$, put $r$ in square $i$. The inverse of $\varphi$ is also straightforward: given a Yamanouchi shuffle tableau $S$, if square $i$ in $S$ contains $r$, then add a square with number $i$ to row $r$ in $\varphi^{-1}(S)$. Finally, sort all the rows in increasing order.
    \end{definition}

    For example, the readers may check that $\phi$ sends the two $D$-compatible tableaux to the left of Figure \ref{fig:compTableau} to the two Yamanouchi shuffle tableaux to the left of Figure \ref{fig:yamaTableau}.

    \begin{proposition}
        $\varphi$ is a bijection between $D$-compatible tableaux and Yamanouchi shuffle tableaux on $D$ such that for any $D$-compatible tableau $T$, the shape of $T$ is the same as the content of $\varphi(T)$.
    \end{proposition}

    \begin{proof}[Proof sketch]
        We will check the following: for any $D$-compatible tableau $T$ and any Yamanouchi shuffle tableau $S$:

        \begin{enumerate}
            \item $\varphi(T)$ is a valid shuffle tableau (Lemma \ref{lem:phiT-valid}),
            \item $\varphi(T)$ is a Yamanouchi shuffle tableau (Lemma \ref{lem:phiT-yama}),
            \item $\varphi^{-1}(S)$ is a valid SYT (Lemma \ref{lem:phiS-valid}),
            \item $\varphi^{-1}(S)$ is $D$-compatible (Lemma \ref{lem:phiS-comp}).
        \end{enumerate}
    \end{proof}

    \begin{lemma}\label{lem:phiT-valid}
        For any $D$-compatible tableau $T$, $\varphi(T)$ is a valid shuffle tableau.
    \end{lemma}

    \begin{proof}
        This follows immediately from the compatibility conditions: the rows of $\varphi(T)$ are weakly increasing thanks to condition (1) in Definition \ref{def:compatible}; the columns of $\varphi(T)$ are strictly increasing thanks to condition (2).
    \end{proof}

    \begin{lemma}\label{lem:phiT-yama}
        For any $D$-compatible tableau $T$, $\varphi(T)$ is a Yamanouchi shuffle tableau.
    \end{lemma}

    \begin{proof}
        Suppose for contradiction that square $i$ containing $r$ in $\varphi(T)$ can be lowered by $E_{r-1}$. This means that $r\geq 2$, so in $T$, there is a $j$ directly above $i$. Recall that we say square $(r_1,c_1)$ is right of square $(r_2,c_2)$, and square $(r_2,c_2)$ is left of square $(r_1,c_1)$, if $r_2 > r_1$ or $r_2 = r_1$ and $c_2 < c_1$. Because $T$ is a SYT, in $T$, the number of entries smaller than $j$ on row $r-1$ equals that of entries smaller than $i$ on row $r$.

        In $\varphi(T)$, this means that the number of entries $r-1$ right of square $j$ equals that of entries $r$ right of $i$. Note that $j$ contains an additional $r-1$ right of $i$. On the other hand, $E_{r-1}$ lowers square $i$, so the number of non-$(r-1,r)$-overlapped $r-1$'s right of square $i$ is less than or equal to that of non-$(r-1,r)$-overlapped $r$'s.

        If none of $j$ or the $(r-1)$-squares right of $j$ is $(r-1,r)$-overlapped with $i$ or some $r$-square left of $i$, then we still have the number of non-$(r-1,r)$-overlapped $r-1$'s right of or equal to $j$ is more than that of non-$(r-1,r)$-overlapped $r$ right of $i$. This contradicts the previous paragraph.

        Thus, either $j$ or some $(r-1)$-square right of $j$ has to be $(r-1,r)$-overlapped with either $i$ or some $r$-square left of $i$. Since $E_{r-1}$ lowers square $i$, $i$ cannot be $(r-1,r)$-overlapped. Thus, this $(r-1,r)$-overlap pair has to be $(s,t)$ where $s$ is $j$ or right of $j$, and $t$ is left of $i$.

        This means that in $T$, $t$ is strictly East of $i$, but $s$ is weakly West of $j$ and $i$, which contradicts condition (2) in Definition \ref{def:compatible}.
    \end{proof}

    \begin{lemma}\label{lem:phiS-valid}
        For any Yamanouchi shuffle tableau $S$, $\varphi^{-1}(S)$ is a valid SYT.
    \end{lemma}

    \begin{proof}
        Row strictness follows from the construction. We now check column strictness. Suppose for contradiction that in $\varphi^{-1}(S)$, square $(r,c)$ (row $r$, column $c$) has entry $i$ and square $(r+1,c)$ has entry $j$ where $i > j$. Then, square $j$ in $S$ contains $r+1$. Furthermore, since the rows are sorted in increasing order, there are $c-1$ more squares containing $r+1$ right of $j$. Since $i > j$, there are at most $c-1$ squares containing $r$ right of $j$. Thus, there are at least as many squares containing $r+1$ than $r$ right of $j$ in $S$. Hence, we can apply $E_r$ on $S$, so $S$ is not Yamanouchi, contradiction.
    \end{proof}

    \begin{lemma}\label{lem:phiS-comp}
        For any Yamanouchi shuffle tableau $S$, $\varphi^{-1}(S)$ is $D$-compatible.
    \end{lemma}

    \begin{proof}
        Condition (1) in Definition \ref{def:compatible} follows from $S$ is row increasing and $\varphi^{-1}(S)$ is a valid SYT. We now check condition (2).

        Suppose for contradiction that there are squares $i$ and $j$ in $S$ such that $j$ is directly below $i$, but in $\varphi^{-1}(S)$, $j$ is not strictly South or weakly West of $i$. Note that by column strictness of $S$, $j$ is still strictly South of $i$, so $j$ has to be strictly East of $i$. Let $(s,t)$ be the coordinate of $i$, then $(s+2,t)$ is the coordinate of $j$. Let $r$ be the number in $j$, and $r-k$ be the number in $i$.

        First, we claim that $k \leq 2$. If $k > 2$, in $\varphi^{-1}(S)$, there is $j > h_1>h_2> i$ such that $h_1$ and $h_2$ are directly above $j$, and $h_2$ is strictly Southeast of $i$. Then, in $S$, the squares $j$, $h_1$, $h_2$, and $i$ appears from left to right. However, $j$ contains $r$ on row $s+2$, and $i$ contains $r-k$ on row $s$, but $h_1$ and $h_2$ contains $r-1$ and $r-2$ respectively, so they are both on row $s+1$. This is impossible since $h_1$ is left of $h_2$.

        For the case $k = 2$, we make an extra assumption: among all possible pairs of squares $i,j$ containing $r-2$ and $r$ violating condition (2), choose the pair with $i$ largest possible. Let $h$ be the number directly above $j$ in $\varphi^{-1}(S)$. By the same argument above, in $S$, square $h$ contains $r-1$ and is on row $s+1$. Since $h$ is also strictly East of $i$, there are more $(r-1)$-squares right of $h$ than $(r-2)$-squares right of $i$. Note that none of the $(r-2)$-squares right of $i$ can be $(r-2,r-1)$-overlapped with an $(r-1)$-square left of $h$. This is because such $(r-2)$-square has to be on row $s$, but $i$ is not $(r-2,r-1)$-overlapped, so any $(r-2)$-square left of $i$ on the same row cannot be so either. This means that the number of non-$(r-2,r-1)$-overlapped $(r-2)$-squares right of $i$ is less than that of non-$(r-2,r-1)$-overlapped $(r-1)$-squares right of $h$. Thus, for $S$ to be Yamanouchi, the has to be non-$(r-2,r-1)$-overlapped $(r-2)$-squares right of $h$ but left of $i$. Let $i'$ be such square, then $i' > i$ and $i'$ is also on row $s$. Moreover, the square directly above $i'$, call it $j'$, is on row $s+2$, and $j' > j$. Then, $i'$ and $j'$ is another pair violating condition (2). This contradicts the choice of $i,j$ being such pair with $i$ largest possible.

        Finally, for the case $k = 1$, we also make an (opposite) extra assumption: among all possible pairs of squares $i,j$ containing $r-1$ and $r$ violating condition (2), choose the pair with $i$ smallest possible. Let $h$ be the number directly below $i$ in $\varphi^{-1}(S)$, then $h$ contains $r$ and is on row $s+1$. By the choice of $i$, $h$ is not $(r-1,r)$-overlapped, and neither is any $r$-square right of $j$ on row $s+2$. Since $h$ is directly above $i$ in $\varphi^{-1}(S)$, in $S$, there are as many $(r-1)$-squares right of $i$ as $r$-squares right of $h$. Any $(r-1)$-square right of $i$ on row $s$ cannot be $(r-1,r)$-overlapped by the choice of $i$. Hence, there are as many non-$(r-1,r)$-overlapped $(r-1)$-squares right of $i$ as non-$(r-1,r)$-overlapped $r$-squares right of $h$. On the other hand, any $(r-1)$-square left of $i$ but right of $h$, including $i$, is on row $s$ and is $(r-1,r)$-overlapped. Thus, to the right of $h$, there are as many non-$(r-1,r)$-overlapped $(r-1)$-squares as $r$-squares. However, $h$ is itself non-$(r-1,r)$-overlapped, so $S$ is not Yamanouchi. This contradiction completes the proof.
    \end{proof}

\subsection{Peelable tableaux}

    Finally, to obtain Theorem \ref{thm:main-thm}, we make the connection between $D$-compatible tableaux and $D$-peelable tableaux. The connection is simple: a $D$-peelable tableau is an inverse standardization of a $D$-compatible tableau.

    \begin{definition}\label{def:standard}
        Given a SSYT $T$, the \textit{standardization} of $T$ is obtained by first replacing all $1$'s from left to right by $1,\ldots,\gamma_1$, then replacing all (original) $2$'s by $\gamma_1+1,\ldots,\gamma_1 + \gamma_2$, etc..
    
        Given $(\lambda/\mu)$, $(\nu/\rho)$, and a SYT $T$. The \textit{inverse standardization} of $T$ with respect to $(\lambda/\mu)\circledast(\nu/\rho)$ is obtained by first replacing $1,\ldots,\lambda_1 - \mu_1$ by $1$, then replacing $\lambda_1-\mu_1 + 1,\ldots,\lambda_1-\mu_1 + \nu_1-\rho_1$ by $2$, then replacing $\lambda_1-\mu_1 + \nu_1-\rho_1+1,\ldots,\lambda_1-\mu_1 + \nu_1-\rho_1+\lambda_2 - \mu_2$ by 3, etc..
    \end{definition}

    \begin{lemma}
        If $T$ is a $D$-compatible tableau, where $D$ is the shuffle diagram of shape $(\lambda/\mu)\circledast(\nu/\rho)$, then the inverse standardization of $T$ with respect to $(\lambda/\mu)\circledast(\nu/\rho)$ is a $D$-peelable tableau. Conversely, if $T$ is a $D$-peelable tableau, then the standardization of $T$ is a $D$-compatible tableau.
    \end{lemma}

    The proof of this lemma is straightforward and is left to the readers. As a consequence, we can construct a bijection $\psi$ between $D$-peelable tableaux and Yamanouchi shuffle tableaux in the same fashion as Definition \ref{def:yama-comp-bijection}.

    \begin{definition}
        Given a $D$-peel tableau $T$, $\psi(T)$ is constructed as follows: if there are $k$ $i$'s in row $r$ of $T$, then in $\psi(T)$, put $k$ $r$'s in row $i$. Sort all the rows in increasing order. The inverse of $\psi$ is defined similarly.
    \end{definition}

    For example, the readers may check that $\psi$ sends the $D$-peelable tableaux in Figure \ref{fig:peelableTableaux} to the Yamanouchi shuffle tableaux in Figure \ref{fig:yamaTableau}.
    
\subsection{Temperley--Lieb immanants}\label{subsec:TL-imm}

    Now, we make a slight detour to give a better description of the Littlewood--Richardson rule for Temperley--Lieb immanants. First, we review the rule from \cite{nguyen2025temperley}. For the background on Hecke algebra and Temperley--Lieb algebra, we refer the readers to \cite[Section 2.2]{nguyen2025temperley} as well as \cite{haiman1993hecke,rhoades2005temperley}. A fact we would like to note is that the basis elements of the Temperley--Lieb algebra can be indexed by noncrossing matchings of $2n$ points.

    Given partitions $\mu = (\mu_1,\ldots,\mu_n)$ and $\nu = (\nu_1,\ldots,\nu_n)$, the \textit{generalized Jacobi--Trudi matrix} $A_{\mu,\nu}$ is the finite matrix whose $i,j$ entry is $h_{\mu_i - \nu_j}$. By convention, we have $h_0 = 1$ and $h_m = 0$ for all $m < 0$.

    Given $\mu$ and $\nu$, we can define two skew shapes $\mu_R/\nu_R$ and $\mu_B/\nu_B$ where
    \[ \mu_R = (\mu_{2i-1}+i-k)_{i = 1}^k, \]
    \[ \nu_R = (\nu_{2i-1}+i-k)_{i = 1}^k, \]
    \[ \mu_B = \left(\mu_{2i}+i-\left\lfloor\dfrac{n+1}{2}\right\rfloor\right)_{i = 1}^k, \]
    \[ \nu_B = \left(\nu_{2i}+i-\left\lfloor\dfrac{n+1}{2}\right\rfloor\right)_{i = 1}^k, \]
    For a shuffle tableau of shape $(\mu_R/\nu_R)\circledast(\mu_B/\nu_B)$, we obtain the Temperley--Lieb type as follows: for each pair
    \[ \begin{ytableau} 
    \none & j \\ 
    i & \none
    \end{ytableau}, \]
    if $i \leq j$, then draw two lines
    \[ \begin{tikzcd}[sep = tiny]
	{} & j \\
	i & {}
	\arrow[no head, from=2-1, to=1-1]
	\arrow[no head, from=1-2, to=2-2]
    \end{tikzcd}; \]
    if $i > j$, then draw two lines
    \[ \begin{tikzcd}[sep = tiny]
	{} & j \\
	i & {}
	\arrow[no head, from=2-1, to=2-2]
	\arrow[no head, from=1-2, to=1-1]
    \end{tikzcd}. \]
    By convention, the squares North and West of the tableau have value $-\infty$, and those South and East of the tableau have value $\infty$. Finally, we put $L_i$ at the end of the $i$th row and $R_i$ at the beginning of the $i$th row.

    \begin{example}
        Let $\mu = (9,9,7,6,4,2)$ and $\nu = (7,4,2,1,0,0)$, then $\mu_R/\nu_R = (7,6,4) / (5,1,0)$ and $\mu_B/\nu_B = (7,5,2) / (2,0,0)$. A shuffle tableaux of shape $(\mu_R/\nu_R)\circledast(\mu_B/\nu_B)$ can be seen below, together with the lines drawn.
        \[\begin{tikzcd}[sep=tiny]
    	&&&&&&&&&& {R_1} & 1 & {} & 3 & {} & {L_1} \\
    	&&&&& {R_2} & 2 & {} & 2 & {} & 2 & {} & 4 & {} & 4 & {L_2} \\
    	&& {R_3} & 1 & {} & 1 & {} & 1 & {} & 2 & {} & 3 & {L_3} \\
    	& {R_4} & 1 & {} & 2 & {} & 4 & {} & 4 & {} & 4 & {L_4} \\
    	{R_5} & 1 & {} & 2 & {} & 2 & {} & 3 & {L_5} \\
    	{R_6} & {} & 3 & {} & 4 & {L_6}
    	\arrow[no head, from=6-3, to=6-4]
    	\arrow[no head, from=5-4, to=5-3]
    	\arrow[no head, from=6-5, to=6-6]
    	\arrow[no head, from=5-6, to=5-5]
    	\arrow[no head, from=6-4, to=6-5]
    	\arrow[no head, from=5-2, to=4-2]
    	\arrow[no head, from=4-3, to=5-3]
    	\arrow[no head, from=5-4, to=4-4]
    	\arrow[no head, from=4-5, to=5-5]
    	\arrow[no head, from=5-6, to=4-6]
    	\arrow[no head, from=4-7, to=5-7]
    	\arrow[no head, from=4-9, to=5-9]
    	\arrow[no head, from=4-8, to=5-8]
    	\arrow[no head, from=5-7, to=5-8]
    	\arrow[no head, from=4-11, to=4-10]
    	\arrow[no head, from=4-3, to=3-3]
    	\arrow[no head, from=3-4, to=4-4]
    	\arrow[no head, from=5-2, to=6-2]
    	\arrow[no head, from=5-1, to=6-1]
    	\arrow[no head, from=6-2, to=6-3]
    	\arrow[no head, from=4-5, to=4-6]
    	\arrow[no head, from=3-6, to=3-5]
    	\arrow[no head, from=4-7, to=4-8]
    	\arrow[no head, from=3-8, to=3-7]
    	\arrow[no head, from=4-9, to=4-10]
    	\arrow[no head, from=3-10, to=3-9]
    	\arrow[no head, from=4-11, to=4-12]
    	\arrow[no head, from=3-12, to=3-11]
    	\arrow[no head, from=3-6, to=2-6]
    	\arrow[no head, from=3-4, to=3-5]
    	\arrow[no head, from=2-7, to=3-7]
    	\arrow[no head, from=3-8, to=2-8]
    	\arrow[no head, from=2-9, to=3-9]
    	\arrow[no head, from=2-11, to=3-11]
    	\arrow[no head, from=3-12, to=2-12]
    	\arrow[no head, from=2-13, to=3-13]
    	\arrow[no head, from=2-11, to=2-12]
    	\arrow[no head, from=1-12, to=1-11]
    	\arrow[no head, from=3-10, to=2-10]
    	\arrow[no head, from=2-9, to=2-10]
    	\arrow[no head, from=2-7, to=2-8]
    	\arrow[no head, from=2-15, to=1-15]
    	\arrow[no head, from=2-13, to=2-14]
    	\arrow[no head, from=1-14, to=1-13]
    	\arrow[no head, from=1-16, to=2-16]
    	\arrow[no head, from=1-14, to=1-15]
    	\arrow[no head, from=2-14, to=2-15]
    	\arrow[no head, from=1-13, to=1-12]
        \end{tikzcd}\]
        This means that the Temperley--Lieb type is
        \[ \resizebox{!}{0.25\textwidth}{
            \begin{tikzpicture}
                \draw (0,0) node[anchor=center] {$L_1$};
                \draw (0,-1) node[anchor=center] {$L_2$};
                \draw (0,-2) node[anchor=center] {$L_3$};
                \draw (0,-3) node[anchor=center] {$L_4$};
                \draw (0,-4) node[anchor=center] {$L_5$};
                \draw (0,-5) node[anchor=center] {$L_6$};
                \draw (3,0) node[anchor=center] {$R_1$};
                \draw (3,-1) node[anchor=center] {$R_2$};
                \draw (3,-2) node[anchor=center] {$R_3$};
                \draw (3,-3) node[anchor=center] {$R_4$};
                \draw (3,-4) node[anchor=center] {$R_5$};
                \draw (3,-5) node[anchor=center] {$R_6$};
                \filldraw[black] (0.5,0) circle (2pt);
                \filldraw[black] (0.5,-1) circle (2pt);
                \filldraw[black] (0.5,-2) circle (2pt);
                \filldraw[black] (0.5,-3) circle (2pt);
                \filldraw[black] (0.5,-4) circle (2pt);
                \filldraw[black] (0.5,-5) circle (2pt);
                \filldraw[black] (2.5,0) circle (2pt);
                \filldraw[black] (2.5,-1) circle (2pt);
                \filldraw[black] (2.5,-2) circle (2pt);
                \filldraw[black] (2.5,-3) circle (2pt);
                \filldraw[black] (2.5,-4) circle (2pt);
                \filldraw[black] (2.5,-5) circle (2pt);
                \draw[thick] [bend right = 90, looseness=1.25] (0.5, -1) to (0.5, 0);
                \draw[thick] [bend right = 90, looseness=1.25] (1.5, 0) to (1.5, -1);
                \draw[thick] [bend right = 90, looseness=1.25] (0.5, -4) to (0.5, -3);
                \draw[thick] [bend right = 90, looseness=1.25] (1.5, -3) to (1.5, -4);
                \draw[thick] (0.5,-2) to (1.5,-2);
                \draw[thick] (0.5,-5) to (1.5,-5);
                
                \draw[thick] [bend right = 90, looseness=1.25] (1.5, -2) to (1.5, -1);
                \draw[thick] [bend right = 90, looseness=1.25] (2.5, -1) to (2.5, -2);
                \draw[thick] [bend right = 90, looseness=1.25] (1.5, -5) to (1.5, -4);
                \draw[thick] [bend right = 90, looseness=1.25] (2.5, -4) to (2.5, -5);
                \draw[thick] (1.5,0) to (2.5,0);
                \draw[thick] (1.5,-3) to (2.5,-3);
            \end{tikzpicture}
        }. \]
    \end{example}    
    
    The Littlewood--Richardson rule for Temperley--Lieb immanants of $A_{\mu,\nu}$ is as follows.

    \begin{theorem}[{\cite[Theorem 6.2]{nguyen2025temperley}}]\label{thm:LR-rule}
        For any partitions $\mu,\nu$, any Temperley--Lieb type $\tau$, and any partition $\lambda$, the coefficient of the Schur function $s_\lambda$ in $\Imm^{\TL}_\tau(A_{\mu,\nu})$ is the number of Yamanouchi shuffle tableaux of shape $(\mu_R/\nu_R)\circledast(\mu_B/\nu_B)$, Temperley--Lieb type $\tau$, and content $\lambda$.
    \end{theorem}

    In particular, the notion of Yamanouchi shuffle tableaux in \cite{nguyen2025temperley} requires checking that no crystal operator $E_i$ can be applied. This makes generating Yamanouchi shuffle tableaux and computing Littlewood--Richardson coefficients time-consuming. With Theorem \ref{thm:main-thm}, we now have a more efficient method.

    \begin{corollary}\label{cor:new-LR-rule}
        For any partitions $\mu,\nu$, any Temperley--Lieb type $\tau$, and any partition $\lambda$, the coefficient of the Schur function $s_\lambda$ in $\Imm^{\TL}_\tau(A_{\mu,\nu})$ is the number of $(\mu_R/\nu_R)\circledast(\mu_B/\nu_B)$-peelable tableaux $T$ of shape $\lambda$ such that $\psi(T)$ has Temperley--Lieb type $\tau$.
    \end{corollary}

\section{Bender--Knuth involutions and symmetry}\label{sec:BK}

    Bender--Knuth (BK) involutions were originally defined on column-strict tableaux by Bender and Knuth \cite{benderknuth1972} and further studied by several authors, e.g., \cite{SCHUTZENBERGER1972,berenstein-zelevenski,berenstein-kirilov,CGP,chiang2024bender,nguyen2023cactus}. 

    \begin{definition}\label{defn:BK_CST}
    The {\it Bender--Knuth involution} $\BK_i$ for $1 \le i \le n-1$ is an involution of the set of SSYT of shape $\lambda$ that sends a SSYT $T$ to the SSYT obtained from the following procedure:
    \begin{enumerate}
        \item Let $S$ be the skew tableau obtained by taking only the squares of $T$ with entry equal $i$ and $i+1$;
        \item Each row of $S$ contains
            \begin{enumerate}
                \item $a$ entries $i$ that are directly above an $i+1$,
                \item $b$ entries $i$ that are alone in their columns,
                \item $c$ entries $i+1$ that are alone in their columns, and
                \item $d$ entries $i+1$ that are directly below an $i$
            \end{enumerate}
        for some $a,b,c,d \geq 0$;
        \item Construct a skew tableau $S'$ by swapping $b$ and $c$ in each row of $S$;
        \item Define $\BK_i(T)$ to be the tableau obtained by replacing $S$ with $S'$ in $T$.
    \end{enumerate}
    If $\alpha = (\alpha_1,\ldots,\alpha_i,\alpha_{i+1},\ldots,\alpha_\ell)$ is the content of $T$, then $\alpha = (\alpha_1,\ldots,\alpha_{i+1},\alpha_i,\ldots,\alpha_\ell)$ is the content of $\BK_i(T)$.
    \end{definition}
    
    \begin{example}
        Consider the action of $\BK_2$ on
        \[ T = \;
        \begin{ytableau} 
            1 & 1 & 1 & 1 & 2 & 2 & 2 & 2 & 3 \\ 
            2 & 2 & 3 & 3 & 3 & 3 \\
            3 & 4 & 4 & 5
            \end{ytableau} \; .
        \]
        The skew tableau $S$ containing only entries $2$ and $3$ in $T$ is
        \[ S = \;
        \begin{ytableau} 
            \none & \none & \none & \none & 2 & 2 & *(yellow) 2 & *(yellow) 2 & *(pink) 3 \\
            2 & *(yellow) 2 & *(pink) 3 & *(pink) 3 & 3 & 3 \\
            3
            \end{ytableau} \; .
        \]
        The new skew tableau $S'$ is
        \[ S' = \;
        \begin{ytableau} 
            \none & \none & \none & \none & 2 & 2 & *(pink) 2 & *(yellow) 3 & *(yellow) 3 \\
            2 & *(pink) 2 & *(pink) 2 & *(yellow) 3 & 3 & 3 \\
            3
            \end{ytableau} \; ,
        \]
        and hence
        \[ \BK_2(T) = \;
        \begin{ytableau} 
            1 & 1 & 1 & 1 & 2 & 2 & 2 & 3 & 3 \\ 
            2 & 2 & 2 & 3 & 3 & 3 \\
            3 & 4 & 4 & 5
            \end{ytableau} \; .
        \]
    \end{example}

    Observe that if $T$ is a $(\lambda/\mu)\circledast(\nu/\rho)$-peelable tableau, then $T$ has content $(\lambda_1-\mu_1,\nu_1-\rho_1,\lambda_2-\mu_2,\nu_2-\rho_2,\ldots)$. Thus, $\BK_1\circ\BK_3\circ\cdots\BK_{2\ell-1}(T)$ has content $(\nu_1-\rho_1,\lambda_1-\mu_1,\nu_2-\rho_2,\lambda_2-\mu_2,\ldots)$, so one may hope that this image is $(\nu/\rho)\circledast(\lambda/\mu)$-peelable. This is indeed true.

    \begin{repthm}{\ref{thm:BK-symm}}
        Let $T$ be a $(\lambda/\mu)\circledast(\nu/\rho)$-peelable tableau, and let $\ell = \ell(\lambda)$. Let
        \[ T' = \BK_1\circ\BK_3\circ\cdots\circ\BK_{2\ell - 1} (T), \]
        then $T'$ is a $(\nu/\rho)\circledast(\lambda/\mu)$-peelable tableau. Note that $T$ and $T'$ have the same shape.
    \end{repthm}

    \begin{proof}
        For a column $c$ and a SSYT $T$, we denote by $T_{c\rightarrow}$ the SSYT formed by columns $c,c+1,\ldots$ of $T$. Suppose $T'$ is not $(\nu/\rho)\circledast(\lambda/\mu)$-peelable. This means that there is a $k$ such that the number of NE matchings of $k$ and $k+2$ in $T'$ is less than required.

        \underline{Case 1}: $k = 2i-1$. Let $a = \rho_i - \rho_{i+1}$, then there is a column $c$ such that in $T'_{c\rightarrow}$, there are at least $a+1$ more $2i+1$ than $2i-1$. Choose the smallest such $c$, then column $c$ contains only $2i+1$ and not $2i-1$ (otherwise there is a smaller $c$). Let $(r,c)$ be the coordinate of this $(2i+1)$-square, and let $(r,c')$ ($c'\leq c$) be the coordinate of the leftmost $(2i+1)$-square on row $r$. Then column $c'$ also only contains $2i+1$ and not $2i-1$. Then, the number of $2i+2$ and $2i$ in $T_{c'\rightarrow}$ is the same as that of $2i+1$ and $2i-1$ in $T'_{c'\rightarrow}$, respectively. Thus, in $T_{c'\rightarrow}$ there are also at least $a+1$ more $2i+2$ than $2i$, contradicting that $T$ is $(\lambda/\mu)\circledast(\nu/\rho)$-peelable.

        \underline{Case 2}: $k = 2i$. Again, let $a = \mu_i - \mu_{i+1}$, then there is a column $c$ such that in $T'_{c\rightarrow}$, there are at least $a+1$ more $2i+2$ than $2i$. Choose the smallest such $c$, then column $c$ contains only $2i+2$ and not $2i$. Let $(r,c)$ be the coordinate of this $(2i+2)$-square, and let $(r,c')$ ($c'\leq c$) be the coordinate of the leftmost $(2i+1)$-square on row $r$. Note that $(r-1,c)$ does not contain $2i$, so $(r-1,c')$ can contain at most $2i-1$. If $(r-1,c')$ does not contain $2i-1$, then the argument follows similarly to Case 1. If $(r-1,c')$ contains $2i-1$, then the number of $2i-1$ in $T_{c'\rightarrow}$ is at most the same as that of $2i$ in $T'_{c'\rightarrow}$, so we get the same contradiction.
    \end{proof}

    The BK involutions are clearly shape-preserving, so Theorem \ref{thm:BK-symm} gives a shape-preserving bijection between $(\lambda/\mu)\circledast(\nu/\rho)$-peelable tableaux and $(\nu/\rho)\circledast(\lambda/\mu)$-peelable tableaux.

    \begin{corollary}
        For any $\lambda,\mu,\nu,\rho,\kappa$, we have $c^{\kappa}_{\lambda/\mu, \nu/\rho} = c^{\kappa}_{\nu/\rho, \lambda/\mu}$.
    \end{corollary}

\section{Schur log-concavity}\label{sec:log-concave}

    We have a characterization of the weights $\nu$ in $P_{\lambda,\mu}$: for all $1 \leq i,j\leq n$, $\nu_i - \nu_j$ is weakly between $\lambda_i -\lambda_j$ and $\mu_i - \mu_j$. This characterization gives the following lemma.

    \begin{lemma}\label{lem:parapiped-condition}
        Suppose $\lambda = (a^k,1^{n-k-1},0)$, and $\nu = \lambda - (e_m + e_{m+1}+\ldots+e_k)$, and $\lambda,\mu,\nu,\rho$ satisfy the condition of Conjecture \ref{conj:log-concavity}. Then,

        \begin{enumerate}
            \item[(i)] $a = \lambda_m > \mu_m \geq \mu_{m+1} \geq \cdots \geq \mu_{k}$,
            \item[(ii)] $\mu_{m-1} > \mu_m$,
            \item[(iii)] $\mu_m - \mu_{\ell-1} \geq \cdots \geq \mu_k - \mu_{\ell-1} \leq a-1$.
        \end{enumerate}
    \end{lemma}

    The next lemma is straightforward from the definition.

    \begin{lemma}\label{lem:all-col}
        Let $T$ be a $\lambda\circledast\mu$-peelable tableau, where $\lambda = (a^k,1^{n-k-1},0)$, then the leftmost $a$ columns of $T$ all contain $1,3,\ldots,2k-1$.
    \end{lemma}

    We say a $(2i-1)$-square of a SSYT $T$ is changeable if it can be increased to $2i$. This means that it is the rightmost $(2i-1)$-square on its row and the square directly below it does not contain $2i$. Similarly, a $2i$-square is changeable if it can be decreased to $2i-1$. We are ready to state our injection.

    \begin{definition}[Injection $\theta$]\label{def:theta}
        Let $T$ be a $\lambda\circledast\mu$-peelable tableau, $\theta(T)$ is obtained from the following procedure.

        \begin{enumerate}
            \item Locate the leftmost changeable $(2m-1)$-square, call this square $s_m$ and change this $2m-1$ to $2m$.
            \item Locate the leftmost changeable $(2m+1)$-square that is weakly East of $s_m$, call this square $s_{m+1}$ and change this $2m+1$ to $2m+2$, repear for $2m+3,\ldots,2k-1$.
            \item If the result of Step 2 is a $\nu\circledast\rho$-peelable tableau, stop. Else, this means that there is no NE matching of $2k-1$ and $2k+1$. Let $s_{k+1}$ be the (unique) $(2k+1)$-square. Locate the rightmost changeable $(2k+2)$-square and call this square $t_{k+1}$. Increase $s_{k+1}$ to $2k+2$ and decrease $t_{k+1}$ to $2k+1$.
            \item If the result of Step 3 is a $\nu\circledast\rho$-peelable tableau, stop. Else, repeat the process with $2k+3$ and $2k+4$, $2k+5$ and $2k+6$, etc. until the result is a $\nu\circledast\rho$-peelable tableau. Call the involved squares $s_{k+2},t_{k+2},s_{k+3},t_{k+2},\ldots$. The step have to stop eventually since $\ell(\lambda) < \infty$.
        \end{enumerate}
    \end{definition}

    \begin{example}\label{ex:ez-theta-ex}
        Let $\lambda = (7,7,7,7,1,1,1,0)$, $\mu = (9,6,6,5,4,3,3,0)$, $\nu = \lambda - (e_2+e_3+e_4) = (7,6,6,6,1,1,1,0)$, and $\rho = (9,7,7,6,4,3,3,0)$ as in Example \ref{ex:strip}. In the following $\lambda\circledast\mu$-peelable tableau, the leftmost changeable $3$ is colored yellow. The leftmost changeable $5$ weakly East of the yellow $3$ is colored yellow, and so is the next changeable $7$:
        \[ \;
        \begin{ytableau} 
            1 & 1 & 1 & 1 & 1 & 1 & 1 & 2 & 2 & 2 & 2 & 2 & 2 \\
            2 & 2 & 2 & 3 & 3 & 3 & 3 & 4 & 4 & 4 \\
            3 & 3 & *(yellow) 3 & 4 & 5 & 5 & *(yellow) 5 & 6 & 6 \\
            4 & 4 & 5 & 5 & 7 & 7 & *(yellow) 7 & 8 \\
            5 & 5 & 6 & 6 & 8 & 10 \\
            6 & 6 & 7 & 7 & 9 \\
            7 & 7 & 8 & 8 & 10 \\
            8 & 10 & 10 & 11 & 12 \\
            12 & 12 & 13 & 14 \\
            14 & 14
            \end{ytableau} \;.
        \]
        Thus, Step 1 and 2 changes the tableau to
        \[ \;
        \begin{ytableau} 
            1 & 1 & 1 & 1 & 1 & 1 & 1 & 2 & 2 & 2 & 2 & 2 & 2 \\
            2 & 2 & 2 & 3 & 3 & 3 & 3 & 4 & 4 & 4 \\
            3 & 3 & *(yellow) 4 & 4 & 5 & 5 & *(yellow) 6 & 6 & 6 \\
            4 & 4 & 5 & 5 & 7 & 7 & *(yellow) 8 & 8 \\
            5 & 5 & 6 & 6 & 8 & 10 \\
            6 & 6 & 7 & 7 & 9 \\
            7 & 7 & 8 & 8 & 10 \\
            8 & 10 & 10 & 11 & 12 \\
            12 & 12 & 13 & 14 \\
            14 & 14
            \end{ytableau} \;.
        \]
        This tableau is $\nu\circledast\rho$-peelable, so we stop.
    \end{example}

    \begin{example}\label{ex:hard-theta-ex}
        With the same partitions as Example \ref{ex:ez-theta-ex}, we pick the following $\lambda\circledast\mu$-peelable tableau:
        \[ \;
        \begin{ytableau} 
            1 & 1 & 1 & 1 & 1 & 1 & 1 & 2 & 2 & 2 & 2 & 2 \\
            2 & 2 & 2 & 2 & 3 & 3 & 3 & 4 & 4 \\
            3 & 3 & 3 & *(yellow) 3 & 4 & 4 & 5 & 6 \\
            4 & 4 & 5 & 5 & 5 & *(yellow) 5 & 6 \\
            5 & 5 & 6 & 6 & 7 & 7 & *(yellow) 7 \\
            6 & 6 & 7 & 7 & 8 & 8 & 9 \\
            7 & 7 & 8 & 8 & 10 & 10 & 11 \\
            8 & 10 & 12 & 12 & 13 \\
            10 & 12 & 14 \\
            14 & 14
            \end{ytableau} \;.
        \]
        After Step 1 and 2, we have
        \[ \;
        \begin{ytableau} 
            1 & 1 & 1 & 1 & 1 & 1 & 1 & 2 & 2 & 2 & 2 & 2 \\
            2 & 2 & 2 & 2 & 3 & 3 & 3 & 4 & 4 \\
            3 & 3 & 3 & *(yellow) 4 & 4 & 4 & 5 & 6 \\
            4 & 4 & 5 & 5 & 5 & *(yellow) 6 & 6 \\
            5 & 5 & 6 & 6 & 7 & 7 & *(yellow) 8 \\
            6 & 6 & 7 & 7 & 8 & 8 & 9 \\
            7 & 7 & 8 & 8 & 10 & 10 & 11 \\
            8 & 10 & 12 & 12 & 13 \\
            10 & 12 & 14 \\
            14 & 14
            \end{ytableau} \;.
        \]
        This is, however, not $\nu\circledast\rho$-peelable, so we proceed to Step 3, which swaps $9$ and the rightmost changeable $10$ West of it:
        \[ \;
        \begin{ytableau} 
            1 & 1 & 1 & 1 & 1 & 1 & 1 & 2 & 2 & 2 & 2 & 2 \\
            2 & 2 & 2 & 2 & 3 & 3 & 3 & 4 & 4 \\
            3 & 3 & 3 & 4 & 4 & 4 & 5 & 6 \\
            4 & 4 & 5 & 5 & 5 & 6 & 6 \\
            5 & 5 & 6 & 6 & 7 & 7 & 8 \\
            6 & 6 & 7 & 7 & 8 & 8 & *(pink) 10 \\
            7 & 7 & 8 & 8 & *(pink) 9 & 10 & 11 \\
            8 & 10 & 12 & 12 & 13 \\
            10 & 12 & 14 \\
            14 & 14
            \end{ytableau} \;.
        \]
        This is still not $\nu\circledast\rho$-peelable, so we repeat Step 4, swapping 11 and 12, then 13 and 14:
        \[ \;
        \begin{ytableau} 
            1 & 1 & 1 & 1 & 1 & 1 & 1 & 2 & 2 & 2 & 2 & 2 \\
            2 & 2 & 2 & 2 & 3 & 3 & 3 & 4 & 4 \\
            3 & 3 & 3 & 4 & 4 & 4 & 5 & 6 \\
            4 & 4 & 5 & 5 & 5 & 6 & 6 \\
            5 & 5 & 6 & 6 & 7 & 7 & 8 \\
            6 & 6 & 7 & 7 & 8 & 8 & 10 \\
            7 & 7 & 8 & 8 & 9 & 10 & *(pink) 12 \\
            8 & 10 & *(pink) 11 & 12 & *(green) 14 \\
            10 & 12 & *(green) 13 \\
            14 & 14
            \end{ytableau} \;.
        \]
        This is finally $\nu\circledast\rho$-peelable.
    \end{example}

    \begin{theorem}
        Given $\lambda,\mu,\nu,\rho$ satisfying the conditions of Conjecture \ref{conj:log-concavity}. The map $\theta$ in Definition \ref{def:theta} is well-defined and is an injection from $\lambda\circledast\mu$-peelable tableaux to $\nu\circledast\rho$-peelable tableaux. In particular, this implies Theorem \ref{thm:log-concavity}.
    \end{theorem}

    \begin{proof}[Proof sketch]
        We will check the following:
        \begin{enumerate}
            \item we can always find the required squares $s_m,s_{m+1},\ldots,s_{k}$ in Step (1) and (2) (Lemma \ref{lem:step12-defined}),
            \item after Step (1) and (2), every $i+2$ is NE matched with some $i$, for $1\leq i \leq 2k-2$ (Lemma \ref{lem:step12-good}),
            \item if Step (3) and (4) are needed, we can always find the required squares $t_{k+1},t_{k+2},\ldots$ (Lemma \ref{lem:step34}),
            \item after Step (3) and (4), the result is $\nu\circledast\rho$-peelable (Lemma \ref{lem:step34-good}),
            \item $\theta$ is an injection (Lemma \ref{lem:injection}).
        \end{enumerate}
    \end{proof}

    The main strategy of proving the first point is as follows. Suppose at some point we cannot find the required $s_{m+i}$, then we want to show that this implies $\mu_m \geq \lambda_m$, which contradicts Lemma \ref{lem:parapiped-condition}. To do that, we will find the \textit{path of $\theta$}, denoted $p$, consisting of squares $p_1,p_2,\ldots$ such that the entries in these squares are weakly increasing, range from $2m$ to $2k$, and these squares occupy consecutive columns, starting from column 1, one square for each column. If at some point we cannot find $s_{m+i}$, the path will have length $\lambda_m$, which will imply $\mu_m\geq \lambda_m$ from the NE matching condition of peelable tableaux.

    \begin{definition}[The path of $\theta$]\label{def:path}
        We construct the path $p$ of $\theta$ as follows. First, scan from left to right through every $(2m-1)$-square from column 1 to $s_m$ if exists, else to column $a$ (recall that $\lambda_m = a$), if any $(2m-1)$-square has a $2m$-square directly above, add that $2m$-square to $p$. Then, for $1\leq i\leq k-m$, if we fail to find $s_{m+i-1}$, stop. Else, let $(r_{m+i-1}, c_{m+i-1})$ be the coordinate of $s_{m+i-1}$. Scan from left to right through every $(2(m+i)-1)$-square from row $r_{m+i-1}+1$ to $s_{m+i}$ if exists, else to column $a$, if any $(2(m+i)-1)$-square has a $2(m+i)$-square directly above, add that $2(m+i)$-square to $p$.
    \end{definition}

    \begin{example}
        Recall the peelable tableaux in Example \ref{ex:ez-theta-ex} and \ref{ex:hard-theta-ex}. Their paths of $\theta$ are shown below in pink (with the irrelevant rows excluded).
        \[ \;
        \begin{ytableau} 
            3 & 3 & *(yellow) 3 & 4 & 5 & 5 & *(yellow) 5 & 6 & 6 \\
            *(pink) 4 & *(pink) 4 & 5 & 5 & 7 & 7 & *(yellow) 7 & 8 \\
            5 & 5 & *(pink) 6 & *(pink) 6 & *(pink) 8 & 10
            \end{ytableau} \; \quad\quad\quad \begin{ytableau} 
            3 & 3 & 3 & *(yellow) 3 & 4 & 4 & 5 & 6 \\
            *(pink) 4 & *(pink) 4 & 5 & 5 & 5 & *(yellow) 5 & 6 \\
            5 & 5 & *(pink) 6 & *(pink) 6 & 7 & 7 & *(yellow) 7 \\
            6 & 6 & 7 & 7 & *(pink) 8 & *(pink) 8 & 9 
            \end{ytableau}.
        \]
    \end{example}

    \begin{lemma}
        The squares $p_1,p_2,\ldots$ in the path $p$ of $\theta$ satisfy: (a) their entries are weakly increasing, range from $2m$ to $2k$, and (b) they are in consecutive columns, starting from column 1, one square for each column.
    \end{lemma}

    \begin{proof}
        (a) follows from the construction, we only need to prove (b). First, the $2m$-squares in $p$ are in consecutive columns, starting from column 1, because there is a $(2m-1)$-square in each of the first $a$ columns. Here is an example with $m = 3$
        \[ \begin{ytableau} 
            \none & \none & \none & \none & \none & \none & 5 & 5 & 5 & *(yellow) 5 \\
            \none & \none & \none & 5 & 5 & 5 & *(pink) 6 \\
            5 & 5 & 5 & *(pink) 6 & *(pink) 6 & *(pink) 6 \\
            *(pink) 6 & *(pink) 6 & *(pink) 6 
        \end{ytableau}. \]
        Recursively, suppose we have found $s_{m+i-1}$ at coordinate $(r_{m+i-1}, c_{m+i-1})$ and $p$ has a square in each of the first $h$ columns. Then, we start scanning the $(2(m+i)-1)$-squares precisely from the square at coordinate $(r_{m+i-1}+1, h+1)$, and since there is also a $(2(m+i)-1)$-square in each of the first $a$ columns, we add the $2m$-squares in consecutive columns starting from $h+1$. Here is an example with $m+i = 5$
        \[ \begin{ytableau}
            \none & \none & \none & \none & \none & \none & 9 & 9 & 9 & *(yellow) 9 \\
            7 & 7 & 7 & *(yellow) 7 & 9 & 9 & *(pink) 10 & *(pink) 10 \\
            *(pink) 8 & 9 & 9 & 9 & *(pink) 10 & *(pink) 10 \\
            9 & *(pink) 10 & *(pink) 10 & *(pink) 10
        \end{ytableau}. \]
    \end{proof}

    \begin{lemma}\label{lem:step12-defined}
        In the process of $\theta$, we can always find the required squares $s_m,s_{m+1},\ldots,s_{k}$ in Step (1) and (2).
    \end{lemma}

    \begin{proof}
        If at some point we cannot find $s_{m+i}$, the path of $\theta$ forms a sequence of $\lambda_m$ increasing entries, each strictly East of the previous. Hence, by the NE matching condition of peelable tableaux, $\mu_m$ is at least $\lambda_m$, contradicting Lemma \ref{lem:parapiped-condition}.
    \end{proof}
    
    \begin{lemma}\label{lem:step12-good}
        After Step (1) and (2), every $i+2$ is NE matched with some $i$, for $1\leq i \leq 2k-2$.
    \end{lemma}

    \begin{proof}
        The case for odd $i$ is straightforward because of Lemma \ref{lem:all-col} and the choice of changeable squares weakly from West to East.

        Let $T'$ be the tableau after Step (1) and (2). If the NE matching condition is not satisfied for some $2i$ and $2i+2$. This means that there is a column $h$ such that there are more $2i+2$ in $T'_{h\rightarrow}$ than $2i$, choose the largest such $h$. Let $(r_i,c_i)$ (resp. $(r_{i_1},c_{i+1})$) be the coordinate of $s_i$ (resp. $s_{i+1}$). Note that for $j \leq c_i$ and $j > c_{i+1}$, the difference of $2i$-entries and $(2i+2)$-entries in $T'_{j\rightarrow}$ is the same as in $T_{j\rightarrow}$. Hence, $c_{i} < h \leq c_{i+1}$, and column $h$ has a $(2i+2)$-square. However, by the choice of $s_{i+1}$, every column between $c_i$ and $h$ (inclusive) also has a $(2i+2)$-square, so in $T_{c_i\rightarrow}$, there has to be more $2i+2$ than $2i$, contradicting that $T$ is peelable.
    \end{proof}

    \begin{lemma}\label{lem:step34}
        Let $T'$ be the tableau after Step (2), if $T'$ is not $\nu\circledast\rho$-peelable, then we can find the required $t_{k+1}, t_{k+2},\ldots$ squares in Step (3) and (4).
    \end{lemma}

    \begin{proof}
        Since there is only one $(2k+1)$-square, namely $s_{k+1}$, for $T'$ to be not $\nu\circledast\rho$-peelable, we need $s_{k+1}$ to be in column $a = \lambda_k$, and the changeable $(2k-1)$-square $s_k$ is also in this column. This means that the path of $\theta$ has a square in each of the first $a-1$ columns. This means that there are at least $a-1$ $2m$-squares. Combining with Lemma \ref{lem:parapiped-condition} (i), we must have $\mu_m = a-1$, and so by (iii), $\mu_{k+1}$ is at least $1$. Thus, we can always find $t_{k+1}$.
        
        In addition, if there is a $2k$-square in column some $c\geq a$, then we can add this square to the path and get a path of length $a$. This leads to the same contradiction as in Lemma \ref{lem:step12-defined}. Hence, any $2k$ is West of column $a$, and so is any $(2k+2)$-square. This means that $s_{k+1}$ is actually East of $t_{k+1}$.

        Finally, suppose we have increased some $s_{k+i}$ and decreased some $t_{k+i}$, yet the tableau is still not $\nu\circledast\rho$-peelable. Suppose also that we have $s_{k+i}$ is East of $t_{k+i}$. We still have $t_{k+i+1}$ exists for the same reason as $t_{k+1}$. We now show that $s_{k+i+1}$ is also East of $t_{k+i+1}$. Let $(r_{k+i},c_{k+i})$ be the coordinate of $t_{k+i}$, and let $(r_{k+i},c'_{k+i})$ be the leftmost $2(k+i)$-square in the same row as $t_{k+i}$. By the choice of $t_{k+i}$ being the rightmost changeable $2(k+i)$-square, there is no $2(k+i)$-square East of column $c'_{k+i}$. Hence, $t_{k+1}$ is weakly West of column $c'_{k+i}$. On the other hand, for the tableau to be not $\nu\circledast\rho$-peelable $s_{k+i+1}$ has to be strictly East of $t_{k+i}$. Then, $s_{k+i+1}$, it also has to be strictly East of $(r_{k+i},c'_{k+i})$, and thus is East of $t_{k+i+1}$.
    \end{proof}

    \begin{lemma}\label{lem:step34-good}
        Let $T'$ be the tableau after Step (4), then $T'$ is $\nu\circledast\rho$-peelable.
    \end{lemma}

    \begin{proof}
        The NE matching condition for $2k,2k+2,2k+4,\ldots$ is straightforward since the new ``even'' squares are at $s_{k},s_{k+1},\ldots$, which are each NE of the next.

        For $2k+1,2k+3,\ldots$, we observe also that $t_{k+1},t_{k+2},\ldots$ are also each NE of the next. If $t_{k+i}$ is not NE of $t_{k+i+1}$ for some $i$, then above $t_{k+i}$ must be a $(2(k+i)+1)$-square, but then we would not need to proceed to $t_{k+i+1}$.
    \end{proof}

    \begin{lemma}\label{lem:injection}
        The map $\theta$ is an injection.
    \end{lemma}

    \begin{proof}
        Every step in $\theta$ is invertible. Suppose $T' = \theta(T)$, then one can recover $T$ uniquely from $T'$ as follows.

        \begin{enumerate}
            \item For $m\leq i \leq k$, in the first $a$ columns, there is a unique column $c_i$ that does not contain $2i-1$. This column contains a $2i$-square. We change this square to $2i-1$.
            \item If the result of Step 1 is not $\lambda\circledast\mu$-peelable, $t_{k+1}$ is the unique $(2k+1)$-square, and $s_{k+1}$ is the rightmost $(2k+2)$-square. Decrease $s_{k+1}$ to $2k+1$ and increase $t_{k+1}$ to $2k+2$.
            \item If the result of Step 2 is a $\lambda\circledast\mu$-peelable tableau, stop. Else, repeat the process with $2k+3$ and $2k+4$, $2k+5$ and $2k+6$, etc. until the result is a $\lambda\circledast\mu$-peelable tableau.
        \end{enumerate}
    \end{proof}

\bibliography{bibliography}
\bibliographystyle{alpha}

\end{document}